\documentclass{ifacconf}

\usepackage[nosectionbib]{apacite}
\usepackage{natbib}
\usepackage{graphicx} 
\usepackage{times} 
\usepackage{amsmath} 
\usepackage{amssymb}  
\usepackage{subfigure}
\usepackage{url}
\usepackage{wrapfig}
\usepackage{caption}
\DeclareCaptionLabelFormat{nonumber}{#1}
\usepackage{bbm, dsfont}
\usepackage{mathtools}
\usepackage{multirow}
\usepackage{centernot}
\usepackage{stmaryrd}
\usepackage{verbatim}
\newtheorem{theorem}{Theorem}

\usepackage{color}
\usepackage{algorithm}
\usepackage{tikz}
\usetikzlibrary{matrix}

\newtheorem{definition}{Definition}
\newtheorem{proposition}{Proposition}
\newtheorem{corollary}{Corollary}
\newtheorem{assumption}{Assumption}
\newtheorem{lemma}{Lemma}
\newtheorem{remark}{Remark}
\newtheorem{example}{Example}

{ 

}

\newcommand{\bi}{\begin{itemize}}
\newcommand{\ei}{\end{itemize}}
\newcommand{\bd}{\begin{displaymath}}
\newcommand{\ed}{\end{displaymath}}
\newcommand{\be}{\begin{eqnarray*}}
\newcommand{\ee}{\end{eqnarray*}}

\newcommand{\K}{{\bf K}}
\newcommand{\U}{{\mathbb{U}}}
\newcommand{\M}{{\mathcal{M}}}
\newcommand{\g}{\gamma}
\usepackage[normalem]{ulem}
\renewcommand*{\thefootnote}{\fnsymbol{footnote}}
\usepackage{ifthen}
\newboolean{showcomments}
\setboolean{showcomments}{true}

\newcommand{\highlight}[1]{\ifthenelse{\boolean{showcomments}} {\textcolor{red}{#1}}{}}
\newcommand{\edited}[1]{\ifthenelse{\boolean{showcomments}} {\textcolor{blue}{#1}}{}}
\newcommand{\sai}[1]{\ifthenelse{\boolean{showcomments}}
{ \textcolor{red}{(Sai says:  #1)}}{}}
\newcommand{\sinha}[1]{\ifthenelse{\boolean{showcomments}}
{ \textcolor{red}{(Subhrajit says:  #1)}}{}}
\newcommand{\enoch}[1]{\ifthenelse{\boolean{showcomments}}
{ \textcolor{red}{(Enoch says:  #1)}}{}}
\begin{document}
\begin{frontmatter}

\title{Koopman Operator Methods for Global Phase Space Exploration of Equivariant Dynamical Systems\thanksref{footnoteinfo}} 

\thanks[footnoteinfo]{This work was supported by a Defense Advanced Research Projects Agency (DARPA) Grant No. DEAC0576RL01830 and an Institute of Collaborative Biotechnologies Grant. The Pacific Northwest National Laboratory (PNNL) is operated by Battelle for the U.S. Department of Energy under Contract DE-AC05-76RL01830.}

\author[First]{Subhrajit Sinha} 
\author[Second]{Sai Pushpak Nandanoori} 
\author[Third]{Enoch Yeung}

\address[First]{Pacific Northwest National Laboratory (e-mail: subhrajit.sinha@pnnl.gov).}
\address[Second]{Pacific Northwest National Laboratory (e-mail: saipushpak.n@pnnl.gov)}
\address[Third]{University of California, Santa Barbara (e-mail: eyeung@ucsb.edu)}

\begin{abstract}                
In this paper, we develop the Koopman operator theory for dynamical systems with symmetry. In particular, we investigate how the Koopman operator and eigenfunctions behave under the action of the symmetry group of the underlying dynamical system. Further, exploring the underlying symmetry, we propose an algorithm to construct a global Koopman operator from local Koopman operators. In particular, we show, by exploiting the symmetry, data from all the invariant sets are not required for constructing the global Koopman operator; that is, local knowledge of the system is enough to infer the global dynamics.
\end{abstract}


\end{frontmatter}

\section{Introduction}\label{section_intro}

Dynamical systems theory is one of the most important branches of mathematics in the sense that it has applications in almost all fields of science and engineering. Any system which changes with time is a dynamical system and hence, they are ubiquitous in nature. As such, both theoretical and numerical analysis of dynamical systems is important. An important class of dynamical systems are the ones which have a symmetry in the sense that there exists some transformations on the state space which carries one solution of the dynamical system to another solution of the dynamical system \cite{field1970equivariant}; \cite{field1980equivariant,golubitsky2012singularities} and the symmetries manifest themselves in asymptotic dynamics, bifurcation, attractor structures etc. \cite{chossat1988symmetry,sparrow2012lorenz,mesbahi_symmetry,koopman_symmetry}. Moreover, symmetries play an important role in synchronization, pattern formation, quantum systems, etc. Mathematically, symmetry is specified by the action of some group on the state space and hence, for studying symmetric dynamical systems, elements of group theory and representation theory are used.

Traditionally, theoretical analysis of dynamical systems is performed by studying the evolution of trajectories in the phase space. However, more recently a different technique is being increasingly used to study dynamical systems, where instead of studying the trajectories in the phase space, the focus is, using transfer operators like Perron-Frobenius operator (P-F) and Koopman operator, on studying the evolution of measures or functions defined on the phase space \cite{Lasota,Vaidya_TAC,Mezic_comparison,mezic2005spectral,Mehta_comparsion_cdc,mezic_koopmanism}.

The main advantage of this approach is the fact that the evolution of measures or functions is linear in the infinite-dimensional space. Moreover, the evolution of functions, which is governed by the Koopman operator, is tailor-made for data-driven analysis of dynamical systems. This is especially useful for analysis of higher dimensional systems like power networks, building systems, biological networks, etc. However, one drawback of using transfer operators is that these are typically infinite-dimensional operators. Hence, for data-driven analysis researchers have developed many different methods for computing the finite-dimensional approximations of these transfer operators and using the developed framework for analysis and control of dynamical systems \cite{Dellnitz_Junge,Mezic2000,froyland_extracting,Junge_Osinga,Mezic_comparison,
Dellnitztransport,mezic2005spectral,Mehta_comparsion_cdc,Vaidya_TAC,
raghunathan2014optimal,susuki2011nonlinear,mezic_koopmanism,
mezic_koopman_stability,surana_observer}; \cite{ yeung2018koopman}; \cite{yeung2017learning, sparse_Koopman_acc,johnson2018class,robust_DMD_ACC}; \cite{robust_DMD_arxiv}; \cite{sinha_online_koopman_arxiv}.

In this paper, we use the Koopman operator framework to study dynamical systems with symmetry. In particular, we analyze some basic properties of symmetric dynamical systems and their symmetry group and investigate how these properties are reflected on the infinite-dimensional Koopman operator for the corresponding symmetric dynamical systems. In particular, we analyze how the symmetry of the underlying system affects the evolution of functions in the function space under the action of the Koopman operator. Moreover, as mentioned before, Koopman operator techniques facilitate the data-driven analysis of dynamical systems. To this end, in this paper, we use the construction technique proposed in \cite{global_koopman_sai_arxiv} to provide a method for constructing the global Koopman operator (defined on the entire phase space) from local Koopman operators (defined on locally invariant sets). In particular, we show that using the symmetry of the underlying system, one does not need to train the local Koopman operators on all the different invariant spaces and hence, one \emph{does not need the data from all the invariant subspaces} for constructing the global operator.

\section{Preliminaries}\label{section_preliminaries}
In this section, we discuss the preliminaries of equivariant systems and transfer operators. 
\subsection{Equivariant Dynamical Systems}
Consider a dynamical system
\begin{equation}\label{sys}
    \dot{x} = f(x)
\end{equation}
where $x\in M\subseteq \mathbb{R}^n$ and $f:M\to M$ is assumed to be at least $\cal C^1$. A symmetry of the dynamical system (\ref{sys}) is a transformation that maps solutions of the system to other solutions of the system. A dynamical system with such a transformation is known as an equivariant dynamical system and is defined as follows:
\begin{definition}[Equivariant Dynamical System]
Consider the dynamical system \eqref{sys} and let $G$ be a group acting on $M$. 
The system is called $G$-equivariant if 
\[f(g\cdot x)=g\cdot f(x), \;\text{for}\; g\in G, x\in M,\]
that is, the following diagram commutes for every $g\in G$.
\end{definition}

\begin{wrapfigure}[8]{l}{0.5\columnwidth}
\vspace{-0.5 cm}
\begin{center}
    \includegraphics[width=0.48\columnwidth]{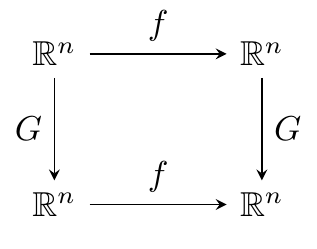}
  \end{center}
  \caption*{}
\end{wrapfigure}
The definition for an equivariant discrete-time dynamical system $x_{t+1}=T(x_t)$ is defined analogously. In particular, a discrete time dynamical system $x_{t+1}=T(x_t)$ is $G$-equivariant if 

\begin{eqnarray}\label{equivariant_discrete}
T(g\cdot x) = g\cdot T(x) \;\text{for}\; g\in G.
\end{eqnarray}
Note that, for a solution $x(t)$ of Eq. \eqref{sys}, $g\cdot x(t)$ is also a solution, with the same being true for discrete-time dynamical system.

\begin{example}
Consider the Lorenz system given by
\begin{eqnarray}\label{sys_lorenz}
\begin{aligned}
 \dot{x} = & \sigma(y-x)\\
 \dot{y} = & x(\rho-z)-y\\
 \dot{z} = & xy-\beta z
\end{aligned}
\end{eqnarray}
where $\sigma$, $\rho$ and $\beta$ are constants. The system equations remain invariant under the transformation $(x,y,z)\mapsto (-x,-y,z)$ and hence the Lorenz system is invariant under the transformation matrix 
\begin{equation}\label{sym_lorenz}
\begin{aligned}
\gamma = \begin{pmatrix}
-1 & 0 & 0\\
0 & -1 & 0\\
0 & 0 & 1\end{pmatrix}.
\end{aligned}
\end{equation}

\begin{figure}[h!]
\begin{center}
\subfigure[]{\includegraphics[width = 0.47 \columnwidth]{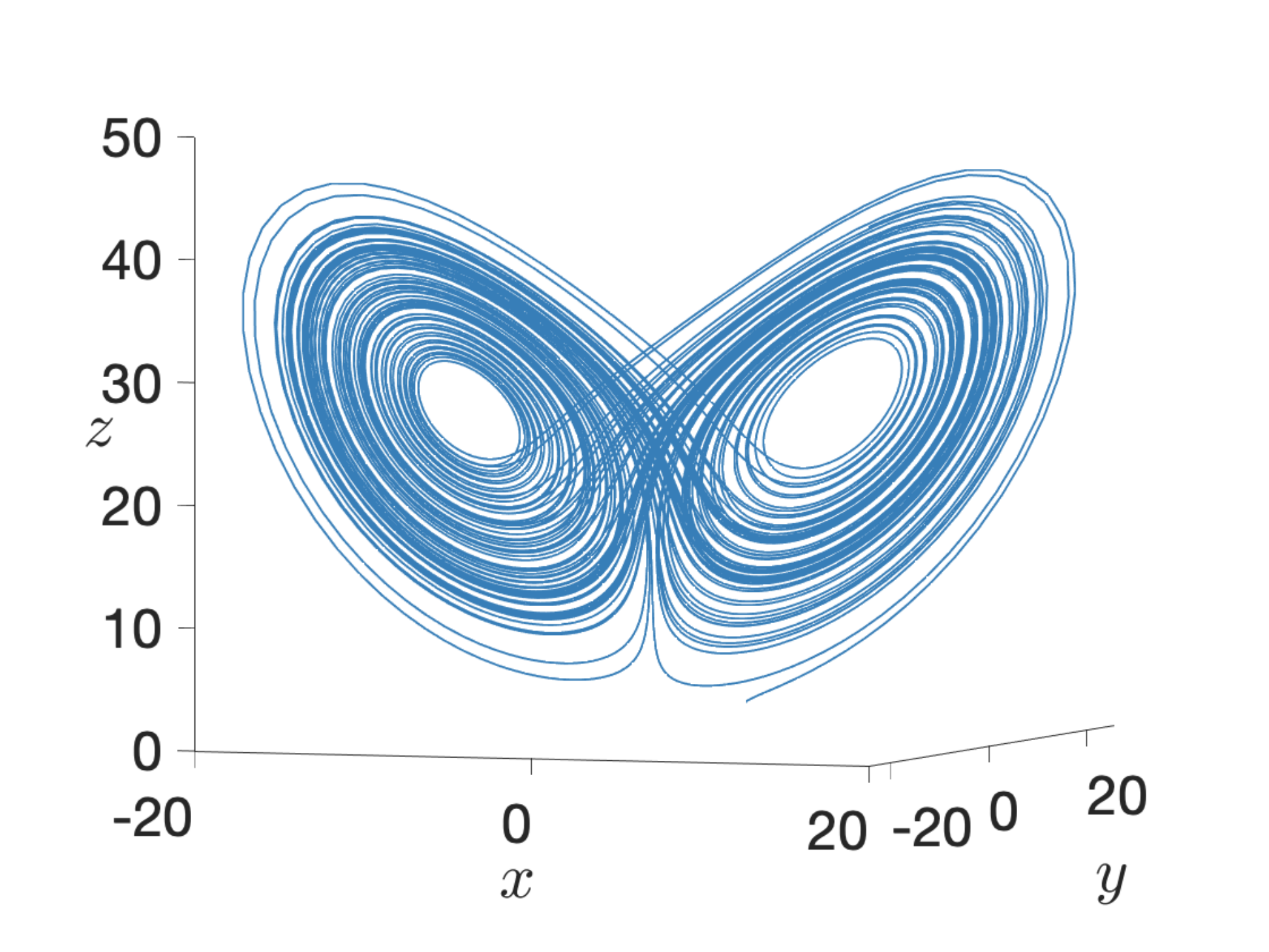}}
\subfigure[]{\includegraphics[width = 0.47 \columnwidth]{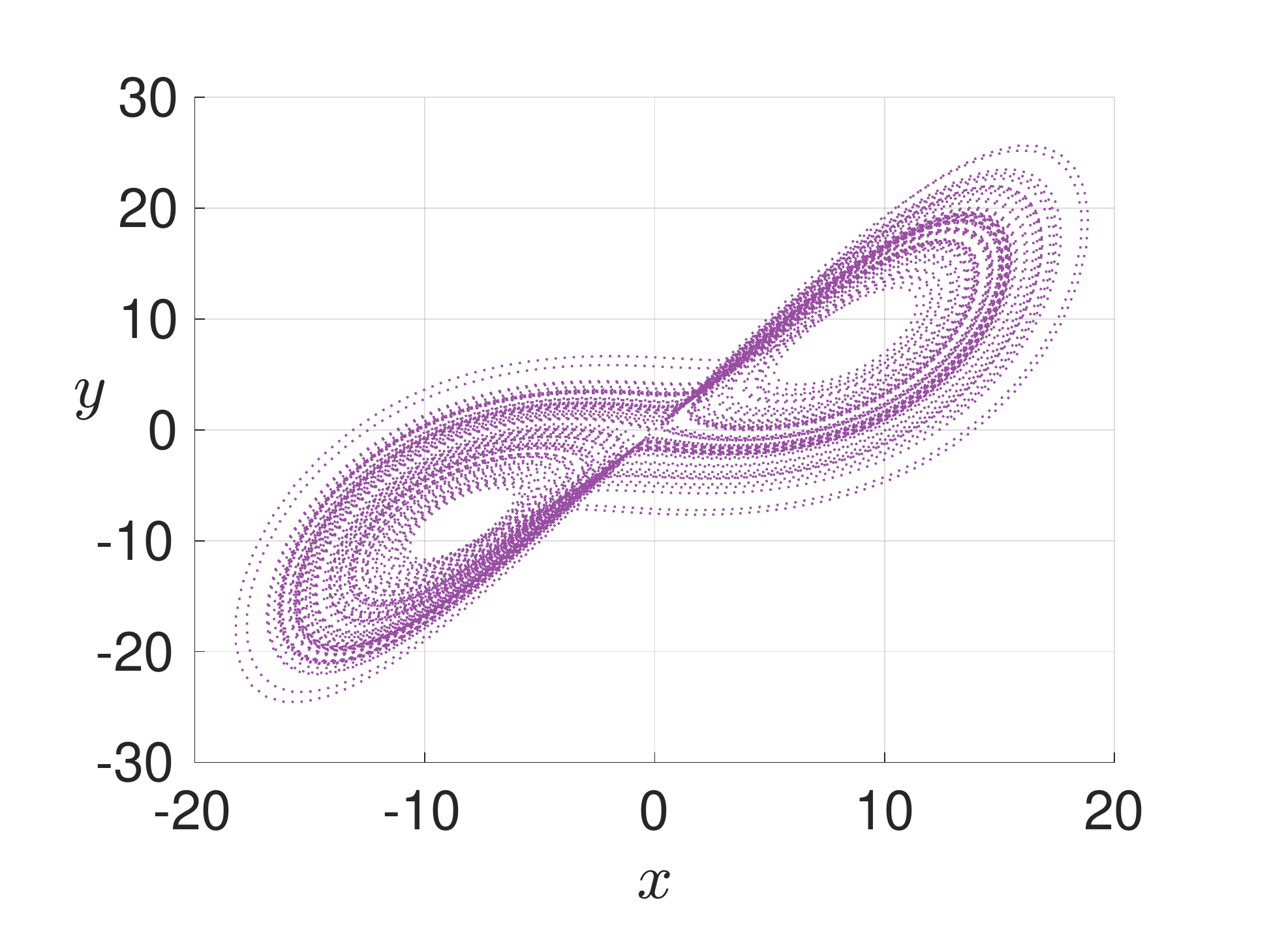}}
\caption{(a) Phase portrait of the Lorenz system (b) Projection of phase portrait on the $x-y$ plane.}
\end{center}
\end{figure}

Note that $\gamma^2 = I$, where $I$ is the identity matrix and this transformation corresponds to a $\pi$ rotation about the $z$-axis. Hence, the Lorenz system is invariant under the group $G = \mathbb{Z}_2$ action, where the action of the non-identity element is given by rotation of $\pi$ about the $z$-axis and $\gamma$ in Eq. (\ref{sym_lorenz}) is the $3$-dimensional representation of the non-identity element of $G$.
\end{example} 
In this paper, we assume $G$ to be a finite subgroup of $O(n)$ and $M\subset\mathbb{R}^n$ to be a compact $G$-invariant set. In general, a symmetry group can be any subgroup of the group of isometries of the Euclidean space $\mathbb{E}^n$, but in this work, we consider finite subgroups of the group of point symmetries of $\mathbb{E}^n$. Further, we consider discrete-time systems of the form
\begin{align}
x_{t+1}=T(x_t).
\label{sys_discrete}
\end{align}
\begin{remark}
We consider discrete-time systems because Koopman operators are tailor-made for data-driven analysis of dynamical systems and data (from a simulation or from an experiment) is always in the form of a discrete time-series. 
\end{remark}
Also, given an abstract group $G$, let $\Gamma$ be the $n$-dimensional representation of the group $G$ in $\mathbb{R}^n$, such that $g\mapsto \gamma_g$, where $g\in G$ and $\gamma_g\in\Gamma$. Note that, $\gamma_g$ has the matrix representation $\gamma_g\in \mathbb{R}^{n\times n}$ and the action of the abstract group $G$ on the state space $\mathbb{R}^n$ is specified by the action of the representation group $\Gamma$ acting on $\mathbb{R}^n$, where the action is by matrix multiplication\footnote{For notational convenience we will use $g$ throughout the paper. However, it should be kept in mind that the action of $g$ is through appropriate representations of $G$ on the concerned spaces.}.

\begin{definition}[Isotropy Set]
Let $x_0(t)$ be a solution (trajectory) of $G$-equivariant dynamical system \eqref{sys_discrete} from the initial condition $x_0$. Then the isotropy set is defined as
\[\Sigma_{x_0(t)}=\{g\in G | g\cdot x_0(t) = x_0(t)\}.\]
\end{definition}

With this we have the following.
\begin{lemma}
The isotropy set corresponding to a solution $x_0(t)$ is a subgroup of the symmetry group $G$.
\end{lemma}
\begin{pf}
Let $g_1,g_2\in\Sigma_{x_0(t)}$. Then, we have 
\[g_1\cdot [g_2\cdot x_0(t)]=g_1\cdot x_0(t)=x_0(t).\]
Therefore, $g_1\cdot g_2\in \Sigma_{x_0(t)}$ and hence $\Sigma_{x_0(t)}$ is closed. Associativity follows directly. Similarly, the identity element belongs to $\Sigma_{x_0(t)}$. Now, suppose $g\in\Sigma_{x_0(t)}$ and let $g_e$ denote the identity element in $G$. Then, we have 
\[x_0(t) = g_e \cdot x_0(t) = g^{-1}g\cdot x_0(t) = g^{-1}\cdot x_0(t).\]
Therefore, for every $g\in\Sigma_{x_0(t)}$, $g^{-1}\in\Sigma_{x_0(t)}$. Hence $\Sigma_{x_0(t)}$ is a subgroup of $G$.
\end{pf}

\begin{proposition}
Let $\Sigma_{x_0(t)}$ and $\Sigma_{g x_0(t)}$ be the isotropy groups of $x_0(t)$ and $g x_0(t)$ respectively. Then we have
\[\Sigma_{g x_0(t)}=g\Sigma_{x_0(t)}g^{-1}.\]
\end{proposition}
\begin{pf}
Let $\theta\in\Sigma_{x_0(t)}$, that is, $\theta$ fixes $x_0(t)$. Then for $g\in G$,
\[(g \cdot \theta \cdot g^{-1})\cdot (g \cdot x_0(t)) = g\cdot \theta\cdot x_0(t) = g\cdot x_0(t).\]
Therefore, $g \cdot \theta \cdot g^{-1}\in\Sigma_{g x_0(t)}$. Similarly, for $\sigma\in\Sigma_{g x_0(t)}$, $g^{-1} \cdot \sigma \cdot g\in\Sigma_{x_0(t)}$. Hence, we obtain
\[\Sigma_{g x_0(t)}=g\Sigma_{x_0(t)}g^{-1}.\]
\end{pf}

\subsection{Transfer Operators}
In this subsection, we briefly discuss the transfer operators, namely the Perron-Frobenius (P-F) and Koopman operator.
Consider a discrete-time dynamical system
\begin{align}\label{system}
x_{t+1} = T(x_t)
\end{align}
where $T: M\subset \mathbb{R}^N \to M$ is assumed to be at least $\cal C^1$.  Associated with the dynamical system (\ref{system}) is the Borel-$\sigma$ algebra ${\cal B}(M)$ on $M$ and the vector space ${\cal M}(X)$ of bounded complex valued measures on $M$. With this, two linear operators, namely, Perron-Frobenius (P-F) and Koopman operator, can be defined as follows \cite{Lasota}:
\begin{definition}
The P-F operator $\mathbb{P}:{\cal M}(X)\to {\cal M}(X)$ is given by
\[[\mathbb{P}\mu](A)=\int_{{M} }\delta_{T(x)}(A)d\mu(x)=\mu(T^{-1}(A))\]
$\delta_{T(x)}(A)$ is stochastic transition function which measure the probability that point $x$ will reach the set $A$ in one time step under the system mapping $T$. 
\end{definition}

\begin{definition}
Invariant measures are the fixed points of
the P-F operator $\mathbb{P}$ that are also probability measures. Let $\bar \mu$ be the invariant measure then, $\bar \mu$ satisfies
\[\mathbb{P}\bar \mu=\bar \mu.\]
\end{definition}

If the state space $X$ is compact, it is known that the P-F operator admits at least one invariant measure.

\begin{definition} 
Given any $h\in\cal{F}$, $\mathbb{U}:{\cal F}\to {\cal F}$ is defined by
\[[\mathbb{U} h](x)=h(T(x))\]
where $\cal F$ is the space of functions (observables) invariant under the action of the Koopman operator.
\end{definition}

\begin{figure}[htp!]
\centering
\includegraphics[scale=.25]{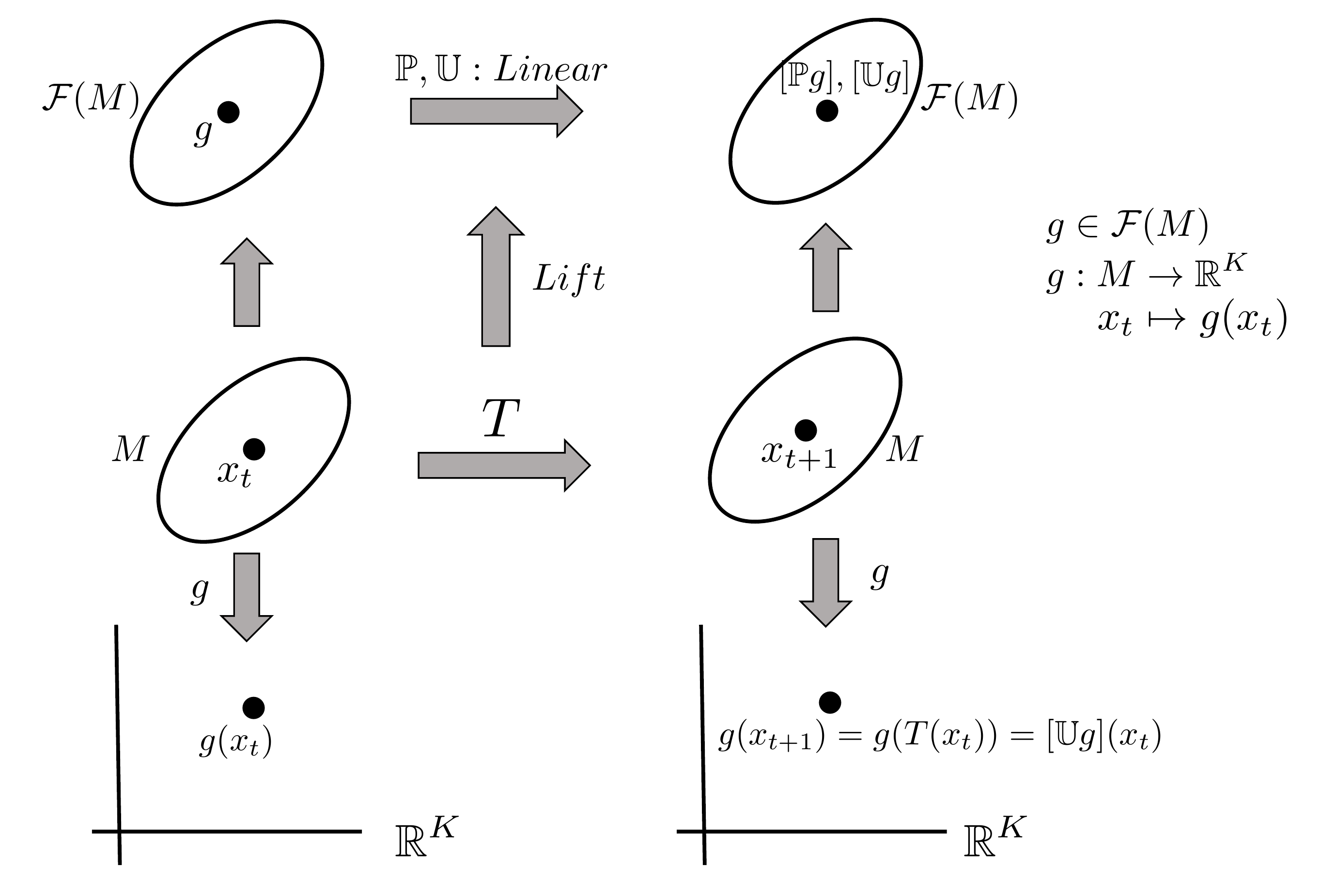}
\caption{Schematic of the P-F and Koopman operators.}\label{koopman_diagram}
\end{figure}

Both the Perron-Frobenius and the Koopman operators are linear operators, even if the underlying system is nonlinear. But while analysis is made tractable by linearity, the trade-off is that these operators are typically infinite-dimensional. In particular, the P-F operator and Koopman operator often will lift a dynamical system from a finite-dimensional space to generate an infinite-dimensional linear system in infinite dimensions (see Fig. \ref{koopman_diagram}).

If the P-F operator is defined to act on the space of densities, that is, $L_1(M)$ and Koopman operator on space of $L_\infty(M)$ functions, then it can be shown that the P-F and Koopman operators are dual to each other\footnote{with some abuse of notation we use the same notation for the P-F operator defined on the space of measure and densities.}.
\begin{align*}
&\left<\mathbb{U} f,g\right>=\int_M [\mathbb{U} f](x)g(x)dx\nonumber\\=&\int_M f(y)g(T^{-1}(y))\left|\frac{dT^{-1}}{dy}\right|dy=\left<f,\mathbb{P} g\right>. 
\end{align*}
where $f\in L_{\infty}(M)$ and $g\in L_1(M)$ and the P-F operator on the space of densities $L_1(M)$ is defined as follows
\[[\mathbb{P}g](x)=g(T^{-1}(x))\left\vert\frac{dT^{-1}(x)}{dx}\right\vert.\]

In the next section, we establish the properties of the Koopman operator for an equivariant dynamical system.
\section{Koopman Operator and Equivariant Dynamical Systems}\label{section_equivariant_koopman}
We begin this section with the analysis of group action on a Koopman operator. 
\subsection{Group Action and Koopman Operator}
Suppose
\begin{align}\label{system_sym}
x_{t+1}=T(x_t)
\end{align}
be a dynamical system defined on the state space $M\subset\mathbb{R}^n$, which is symmetric with respect to a group $G$ and let $\mathbb{U}$ be the associated Koopman operator. The Koopman operator, $\mathbb{U}$ is a linear operator on the space of functions $({\cal F}(M))$ on $M$. We define a map 
\begin{equation}\label{group_action_on_F}
\begin{aligned}
\varphi : & G\times {\cal F}(M)\to {\cal F}(M)\\
& (g\star f)(x) \mapsto f(g^{-1}\cdot x).
\end{aligned}
\end{equation}
\begin{lemma}
The map $\varphi$, defined in Eq. (\ref{group_action_on_F}) defines a group action on the space ${\cal F}(M)$
\end{lemma}
\begin{pf}
Firstly, let $e$ be the identity element of $G$. Then
\[(e\star f)(x)=f(e^{-1}\cdot x)=f(x).\]
Secondly, let $g_1,g_2\in G$ and let $g_1\odot g_2 = g_3\in G$, where $\odot$ denotes the group operation. Then
\begin{eqnarray*}
&&[g_1\star (g_2\star f)](x) = [g_1\star f](g_2^{-1}\cdot x) = f(g_1^{-1}\cdot (g_2^{-1}\cdot x))\\
&&\qquad = f((g_2^{-1}\odot g_1^{-1})\cdot x)=f((g_1\odot g_2)^{-1}\cdot x)\\
&& \qquad = [(g_1\odot g_1)\star f](x).
\end{eqnarray*}
Hence the map $\phi$, defined in (\ref{group_action_on_F}) defines a group action on ${\cal F}(M)$.
\end{pf}


Now a Koopman operator is a linear operator on ${\cal F}(M)$. Hence, from the action of the symmetry group $G$ on ${\cal F}(M)$, we have the following theorem. 

\begin{theorem}
Let $\mathbb{U}$ be the Koopman operator associated with a $G$-equivariant system $x_{t+1}=T(x_t)$. Then
\begin{align}\label{Koopman_commutation}
[g\star (\mathbb{U}f)](x) = [\mathbb{U}(g\star f)](x).
\end{align}
for all $g\in G$ and $f\in{\cal F}(M)$.
\end{theorem}
\begin{pf}
For the dynamical system $x_{t+1}=T(x_t)$ and any function $f\in{\cal F}(M)$, the Koopman operator $\mathbb{U}$ is defined as $[\mathbb{U}f](x)=f(T(x))$. Hence, for $g\in G$ we have,
\begin{align*}
g\star (\U f)(x) & =g \star f(T(x)) = f(g^{-1}\cdot T(x))\\
& =f(T(g^{-1}\cdot x)) = \U f(g^{-1}\cdot x) = [\mathbb{U}(g\star f)](x),
\end{align*}
where the third equality follows from the definition of $G$-equivariant systems (refer Eq. (\ref{equivariant_discrete})).
\end{pf}
The above theorem essentially says that the Koopman operator commutes with the elements of the symmetry group $G$. 

Associated with a Koopman operator is its eigenspectrum, that is, the eigenvalues $\lambda$, and their corresponding eigenfunctions $\phi_{\lambda}(x)$, such that
\[[\U \phi_{\lambda}](x) = \lambda \phi_{\lambda}(x).\]

The eigenspectrum (especially eigenfunctions corresponding to dominant eigenvalues) of a Koopman operator dictates the evolution of the functions $f\in{\cal F}(M)$, under the map $T$ (refer Eq. \ref{system_sym}) and hence the different algorithms like Dynamic Mode Decomposition (DMD) \cite{DMD_schmitt} and Extended Dynamic Mode Decomposition (EDMD) \cite{EDMD_williams} are geared towards obtaining finite-dimensional approximations of the eigenspectrum of the Koopman operator. 

\begin{definition}
Let $\U$ be a Koopman operator and let $\phi_\lambda^i(x)$ be eigenfunctions of $\U$ corresponding to the eigenvalue $\lambda$, that is, $\U\phi_\lambda^i(x)=\lambda\phi_\lambda^i(x)$. Then the eigenspace ${\cal E}_\lambda$ is defined as 
\[{\cal E}_\lambda=\textnormal{span}\{\phi_\lambda^i(x)\}.\]
\end{definition}

The following result establishes that the eigenspace is left invariant under group action. 
\begin{proposition}
Let (\ref{system_sym}) be a $G$-equivariant discrete-time dynamical system and $\U$ be the associated Koopman operator. If $\lambda$ is an eigenvalue of the Koopman operator $\U$ and ${\cal E}_{\lambda}$ is the corresponding eigenspace, then the eigenspace remains invariant under the action of the symmetry group $G$.
\end{proposition}
\begin{pf}
Let $\phi(x)$ be an eigenfunction of the Koopman operator $\U$ with eigenvalue $\lambda$. Then for $g\in G$, we have
\begin{align}
\U (g\star \phi(x)) = g\star \U \phi(x) = g \star \lambda \phi(x) = \lambda (g \star \phi(x))
\end{align}
Hence, if $\phi(x)\in {\cal E}_{\lambda}$, then $(g\star \phi(x))\in {\cal E}_{\lambda}$. Hence the eigenspace is left invariant under the $G$-action.
\end{pf}

Note that a Koopman operator is a linear operator which gives the evolution of functions which are defined on the state space. Let $x\in M$ and $g\cdot x \in M$ for a $G$-equivariant system (\ref{sys_discrete}) and let $f\in L_2(M)$. Let $\hat{f}=g\star f$. Then the following proposition relates the representation (analogous to a matrix representation of a linear transformation) of the Koopman operator when the functions $f$ and $\hat{f}$ are evaluated at $x$ and $g\cdot x$ respectively.
\begin{proposition}\label{psi_i_not_equal_to_psi_j}
Let (\ref{system_sym}) be a $G$-equivariant dynamical system and its associated Koopman operator be $\U:{\cal F}(M)\to {\cal F}(M)$. Suppose $f\in {\cal F}(M)$ and let $\U_f$ be the representation of $\U$ with respect to $f$. For $g\in G$,  let $\hat{f} = g\star f$ and $\U_{\hat{f}}$ be the representation of $\U$ with respect to $\hat{f}$. Then for $x\in M$, we have
\begin{align}\label{koopman_conjugation}
\U_{\hat{f}} \hat{f}(g\cdot x) = \U_f f(x)
\end{align}
\end{proposition}
\begin{pf}
We have
\begin{align*}
\U_{\hat{f}} \hat{f}(g\cdot x) =& \hat{f}(T(g\cdot x)) = \hat{f}(g\cdot T(x))= g^{-1}\star \hat{f}(T(x))\\
=& f(T(x)) = \U_f f(x).
\end{align*}
\end{pf}

\subsection{Group Action and Invariant Spaces}
\begin{definition}
For a dynamical system $x_{t+1}=T(x_t)$, defined on $M\subseteq \mathbb{R}^n$, a subset ${\cal M}\subset M$ is an invariant set if for every trajectory $x_0(t)$, \begin{equation}
x_0(t)\in {\cal M} \implies x_0(\tau) \in {\cal M}, \forall \tau \geq t.
\end{equation}
\end{definition}
Note that an orbit from an initial condition $x_0$ is an invariant set. 

For a measure preserving transformation $T$, all the eigenvalues of the associated Koopman operator ${\U}$ lie on the unit circle \cite{mezic_koopmanism}. Moreover, when $T$ is an ergodic transformation, then all eigenvalues of ${\U}$ are simple \cite{petersen_ergodic_book,mezic_koopmanism}. However, if $T$ is not ergodic, then the state space can be partitioned into subsets ${\cal M}_i$ (minimal invariant subspaces) such that the restriction $T|_{{\cal M}_i}: {\cal M}_i \to {\cal M}_i$ is ergodic. A partition of the state space into invariant sets is called an ergodic partition or stationary partition. Hence, for any transformation $T$, defined on $M\subseteq \mathbb{R}^n$, the state space $M$ can be expressed as 
\begin{align}\label{invariant_union}
M=\cup_{i = 1}^m {\cal M}_i \textnormal{ (modulo measure zero sets)},
\end{align}
where each ${\cal M}_i$ is an invariant set and ${\cal M}_i$ and ${\cal M}_j$ are disjoint for $i\neq j$. Hence, all ergodic partitions are disjoint and they support mutually singular functions from $L_2(M)$ \cite{mezic_koopmanism}. Therefore, the number of linearly independent eigenfunctions of ${\U}$ corresponding to an eigenvalue $\lambda$ is bounded above by the number of ergodic sets in the state space \cite{mezic_koopmanism}. The dynamics of the system dictates the number of ergodic partitions (invariant sets) in the state space. The following results summarize the above discussion. 
\begin{lemma}
Let $\lambda$ be an eigenvalue of the Koopman operator ${\cal K}$.  Suppose if the algebraic multiplicity of the eigenvalue $\lambda$, is equal to the geometric multiplicity, then the corresponding eigenfunctions are linearly independent. 
\label{lemma_am_gm} 
\end{lemma}
\begin{pf}
The proof follows from standard results on matrices \cite{horn2012matrix}. 
\end{pf}

\begin{definition}
Let $\M_i$ be an invariant set of the $G$-equivariant dynamical system (\ref{system_sym}). Then for $g\in G$, define the set $g\cdot\M_i$ as 
\[g\cdot\M_i := \{\Tilde{x}\in M | \Tilde{x} = g\cdot x \textnormal{ for } x\in \M_i\}\]
\end{definition}

\begin{proposition}\label{group_action_invariance}
If $\M_i$ is an invariant set for a $G$-equivariant dynamical system $x_{t+1}=T(x_t)$, then $g\cdot\M_i$ is also an invariant set for $g\in G$.
\end{proposition}
\begin{pf}
For $g\in G$, from the definition of $g\cdot \M_i$, for any $\tilde{x}\in g\cdot\M_i$, there exists some $x\in\M_i$ such that $\tilde{x}=g\cdot x$. Let $\tilde{x}_0\in g\cdot\M_i$ such that $\tilde{x}_0=g\cdot x_0$ for some $x_0\in \M_i$. Since $\M_i$ is an invariant set, the trajectory $x_0(t)\in\M_i$. Now,
\begin{align*}
\tilde{x}_0(t) =& \cup_t \{\tilde{x}\in M|\tilde{x}=T^t(\tilde{x}_0), t\geq 0\}\\
=& \cup_t \{\tilde{x}\in M|\tilde{x}=T^t(g\cdot x_0), t\geq 0\}\\
=& \cup_t \{\tilde{x}\in M|\tilde{x}=g\cdot T^t(x_0), t\geq 0\}
\end{align*}
Now, since $\M_i$ is an invariant set, $T^t(x_0)\in\M_i$ and hence, $\tilde{x}_0(t)\in g\cdot\M_i$ and therefore $g\cdot\M_i$ is invariant. 
\end{pf}
\begin{corollary}
For any invariant set $\M_i$, $G\cdot\M_i$ is invariant, where 
\[G\cdot\M_i =\{\tilde{x}\in M|\tilde{x}=g\cdot x, \textnormal{ for } x \in \M_i \textnormal{ and } g\in G\}.\]
\end{corollary}

\section{Global Phase Space Reconstruction from Data}\label{section_global_phase_space}

In this section, we develop the data-driven tools for analysis of equivariant dynamical systems. 

\subsection{Finite Dimensional Approximation of Koopman Operator}\label{subsection_EDMD}
Let
\begin{align}
X_p = [x_1,x_2,\ldots,x_M], \quad  X_f = [y_1,y_2,\ldots,y_M] \label{data}
\end{align}
be snapshots of data obtained from simulating a dynamical system $x\mapsto T(x)$, or from an experiment, where $x_i\in M$ and $y_i\in M$, $M\subset \mathbb{R}^n$. The two pairs of data sets are assumed to be two consecutive snapshots i.e., $y_i=T(x_i)$. Let $\mathcal{D}=
\{\psi_1,\psi_2,\ldots,\psi_K\}$ be the set of observables, where $\psi_i : M \to \mathbb{R}$ and $\psi_i\in L_2(M)$. Let ${\cal G}_{\cal D}$ denote the span of ${\cal D}$. Let $\mathbf{\Psi}:X\to \mathbb{R}^{K}$ be a vector valued function, such that
\begin{equation*}
\mathbf{\Psi}(\boldsymbol{x}):=\begin{bmatrix}\psi_1(x) & \psi_2(x) & \cdots & \psi_K(x)\end{bmatrix}^\top.
\end{equation*}
Here $\mathbf{\Psi}$ is the mapping from physical space to feature space. Any function $\phi,\hat{\phi}\in \mathcal{G}_{\cal D}$ can be written as
\begin{align}
\phi = \sum_{k=1}^K a_k\psi_k=\boldsymbol{\Psi^T a},\quad \hat{\phi} = \sum_{k=1}^K \hat{a}_k\psi_k=\boldsymbol{\Psi^T \hat{a}}
\end{align}
for some set of coefficients $\boldsymbol{a},\boldsymbol{\hat{a}}\in \mathbb{R}^K$. Let $\hat{\phi}(x)=[\mathbb{U}\phi](x)+r,$
where $r$ is a residual that appears because $\mathcal{G}_{\cal D}$ is not necessarily invariant to the action of the Koopman operator. The finite dimensional approximate Koopman operator $\bf K$ minimizes this residual $r$ and the matrix $\bf K$ is obtained as a solution of the following least square problem: 
\begin{equation}\label{edmd_op}
\min\limits_{\bf K}\parallel{\bf K} {Y_p}-{Y_f}\parallel_F
\end{equation}
where
\begin{align}\label{edmd1}
 {Y_p}={\bf \Psi}(X_p) = & [{\bf \Psi}(x_1), {\bf \Psi}(x_2), \cdots , {\bf \Psi}(x_M)]\\
 {Y_f}={\bf \Psi}(X_f) = & [{\bf \Psi}(y_1), {\bf \Psi}(y_2), \cdots , {\bf \Psi}(y_M)],
\end{align}
with ${\bf K}\in\mathbb{R}^{K\times K}$. The optimization problem (\ref{edmd_op}) can be solved explicitly to obtain following solution for the matrix $\bf K$
\begin{align}
{\bf K}={Y_f}{Y_p}^\dagger \label{EDMD_formula}
\end{align}
where ${Y_p}^{\dagger}$ is the pseudo-inverse of matrix $Y_p$.
DMD is a special case of EDMD algorithm with ${\bf \Psi}(x) = x$.

\subsection{Global Koopman Operator from Local Koopman Operators}\label{subsection_global_koopman_construction}
As mentioned earlier, any phase space $M$ can be decomposed into disjoint invariant sets $\M_i$. Let
\begin{align}\label{dictionary_Mi}
{\bf \Psi}_i = \{\psi_1^i,\psi_2^i,\cdots ,\psi_{K_i}^i\}
\end{align}
be the dictionary functions (observables) on each $\M_i$ and let $\K_i$ be the corresponding Koopman operator on $\M_i$. Note that, in general, on each ${\cal M}_i$, the dictionary functions are different and hence, the finite-dimensional matrix representation of each of the local Koopman operators $\K_i$ is different. This is because given any linear transformation ${\cal O}:V\to W$, where $V$ and $W$ are vector spaces, the matrix representation of $\cal O$ depends on the choice of the basis vectors of $V$ and $W$. For the self-containment of the paper, in this subsection, we briefly review the results of \cite{global_koopman_sai_arxiv} where we had proposed a systematic method to construct the global Koopman operator $\K$, which describes the evolution of the system on the entire phase space $M$, from the local Koopman operators $\K_i$. 

Let ${\bf \Psi}_i$ be the dictionary functions on each invariant set $\M_i$, $i = 1,2,\cdots , m$ and let $\K_i$ be the corresponding Koopman operator which describes the evolution of the system in each ${\cal M}_i$. We define the set of dictionary functions on the entire state space $M$ as
\[{\bf \Psi} = \bigsqcup_{i = 1}^m {\bf \Psi}_i,\]
where $\bigsqcup$ is the disjoint union.
Then if $\K$ is the global Koopman operator on the entire state space $M$ with dictionary function $\bf \Psi$, then $\K$ can be expressed as
\begin{eqnarray}\label{global_Koopman}
\K = \textnormal{diag}(\K_1,\K_2,\cdots , \K_m).
\end{eqnarray}

\subsection{Global Koopman Operator for Equivariant Systems}\label{subsection_global_koopman_equivariant}

Consider the $G$-equivariant system (\ref{system_sym}), defined on the state space $M\in\mathbb{R}^n$ with disjoint invariant sets $\M_i$, as in Eq. (\ref{invariant_union}). 

From proposition \ref{group_action_invariance} we have that $g\cdot \M_i$ is also invariant for all $g\in G$. Hence, $g\cdot \M_i\subset\M_j$ for some $j\in\{1,2,\cdots , m\}$. 
\begin{assumption}
We assume there exists some $g\in G$, such that $g\cdot \M_i\subset\M_j$ and $i\neq j$.
\end{assumption}

Let ${\bf \Psi}_i$ be the dictionary functions defined on $\M_i$ and let $\K_i$ be the local Koopman operator on $\M_i$. Let
\[X=[x_1,x_2,\cdots , x_{M+1}]\]
be points in $\M_i$, such that $x_{k+1}=T(x_k)$ and thus 
\begin{align}\label{koopman_Ki}
[\K_i {\bf \Psi}_i](x_k)={\bf \Psi}_i(x_{k+1}).
\end{align}
Now, $\K_j$ governs the evolution of dictionary functions on $\M_j$ and the goal is to compute $\K_j$. Since in computing $\K_i$, ${\bf \Psi}_i(x_k)$ are already computed, we would like to use this information for computation of $\K_j$. This can be done in two different ways.

\emph{Case I.} We use the same dictionary function ${\bf \Psi}_i$ on $\M_j$, that is 
\[{\bf \Psi}_j={\bf \Psi}_i.\]
\begin{theorem}\label{K_i_K_j_theorem}
Let $\M_i$ be an invariant set of the $G$-equivariant system (\ref{system_sym}) and let $\K_i\in\mathbb{R}^{K_i\times K_i}$ be the local Koopman operator on $\M_i$ with dictionary function ${\bf \Psi}_i(x)$, $x\in\M_i$. Let, for $g\in G$, $g\cdot\M_i\subset \M_j$. Let $\K_j$ be the local Koopman operator on $\M_j$ with dictionary functions ${\bf \Psi}_j={\bf \Psi}_i$. Then 
\[\K_i = \g\K_j \g^{-1},\]
where $g\mapsto \g\in \Gamma$ and $\Gamma$ is the $K_i$ dimensional matrix representation of $G$ in $\mathbb{R}^{K_i}$.
\end{theorem}
\begin{pf}
Since ${\bf \Psi}_j={\bf \Psi}_i$, $\K_j\in\mathbb{R}^{K_i\times K_i}$. Now consider $x_k\in\M_i$ and $g\cdot x_k\in\M_j$. Then we have,
\begin{align}
\K_j{\bf \Psi}_i(g\cdot x_k) =& {\bf \Psi}_i(g\cdot x_{k+1}) = \g^{-1} {\bf \Psi}_i(x_{k+1}) \nonumber \\
=& \g^{-1}\K_i{\bf \Psi}_i(x_k).\label{K_i_K_j1} 
\end{align}
Consider again,
\begin{align}\label{K_i_K_j2}
\K_j{\bf \Psi}_i(g\cdot x_k) = \K_j \g^{-1}{\bf \Psi}_i(x_k).
\end{align}
Hence, from Eqs. (\ref{K_i_K_j1}) and (\ref{K_i_K_j2}), we have
\begin{align}\label{psi_i_equal_to_psi_j}
\K_i{\bf \Psi}_i(x_k) = \g\K_j \g^{-1}{\bf \Psi}_i(x_k).
\end{align}
Since Eq. (\ref{psi_i_equal_to_psi_j}) is true for all $x_k\in\M_i$, we obtain
\[\K_i = \g\K_j \g^{-1}.\]
\end{pf}

\begin{corollary}\label{K_i_K_j_corollary}
Let $x_0(t)$ be a trajectory of a $G$-equivariant system $x_{t+1}=T(x_t)$ and let $g\cdot x_0(t)$ be the image of $x_0(t)$ under the action of $g\in G$. Let ${\bf \Psi}\in L_2(M)$ be a set of dictionary functions of cardinality $K$. Let $\K_{x_0}$ and $\K_{g\cdot x_0}$ be the finite-dimensional representation of Koopman operator which governs the evolution of ${\bf \Psi}$ on $x_0(t)$ and $g\cdot x_0(t)$ respectively. Then
\begin{align}\label{K_i_K_j_trajectory}
\K_{x_0} = \g\K_{g\cdot x_0}\g^{-1},
\end{align}
\end{corollary}
where $g\mapsto \g\in \Gamma$ and $\Gamma$ is the $K$ dimensional matrix representation of $G$ in $\mathbb{R}^{K}$.
\begin{pf}
Similar to the proof in theorem \ref{K_i_K_j_theorem}.
\end{pf}
Note that the local Koopman operators obtained using the DMD algorithm satisfy theorem \ref{K_i_K_j_theorem}.

\emph{Case II.} We define the dictionary function on $g\cdot\M_i$ and hence on $\M_j$, $j\neq i$ as
\[{\bf \Psi}_j = g\star {\bf \Psi}_i.\]
In this case, from proposition \ref{psi_i_not_equal_to_psi_j}, we have \[\K_i=\K_j.\]

Hence, starting with an invariant set $\M_i$, with a local Koopman operator $\K_i$, if $g_j\cdot\M_i\subset\M_j$ for $j\in\{1,2,\cdots , m\}$, $j\neq i$, we can obtain all the local Koopman operators $\K_j$, $j=1,2,\cdots , m$. Hence, using the procedure of \cite{global_koopman_sai_arxiv}, we can obtain the global Koopman operator, defined on the entire state space $M$ for the $G$-equivariant system (\ref{system_sym}).

\section{Simulation Results}\label{section_simulation}
In this section, by applying the symmetry in the system, we identify the global Koopman operator starting from an invariant subspace. We begin with the discussion on systems with reflective symmetry. In all the examples, we use the same dictionary functions on two different invariant spaces (or two different trajectories), which are related by the symmetry group (as in theorem \ref{K_i_K_j_theorem}).

\subsection{Reflection Symmetry: Bistable Toggle Switch}
Consider the bistable toggle switch system, first introduced in \cite{gardner2000construction}. The governing equations of motion are: 
\begin{equation}
\begin{aligned}
\dot{x}_1 = & \frac{\alpha_1}{1+x_2^{\beta}} - \kappa_1 x_1 \\
\dot{x}_2 = & \frac{\alpha_2}{1+x_1^{\theta}} - \kappa_2 x_2
\end{aligned}
\label{eq:bistable_toggle_switch}
\end{equation}
where the states $x_1 \in \mathbb{R}$ and $x_2 \in \mathbb{R}$ indicate the concentration of the repressor $1$ and $2$; the effective rate of synthesis of repressor $1$ and $2$ are denoted by $\alpha_1$ and $\alpha_2$; the self decay rates of concentration of repressor $1$ and $2$ are given by $\kappa_1>0$ and $\kappa_2>0$; the cooperativity of repression of promoter $2$ and $1$ are respectively denoted by $\beta$ and $\theta$. 

\begin{figure}[h!]
\begin{center}
\includegraphics[width = 0.95 \columnwidth]{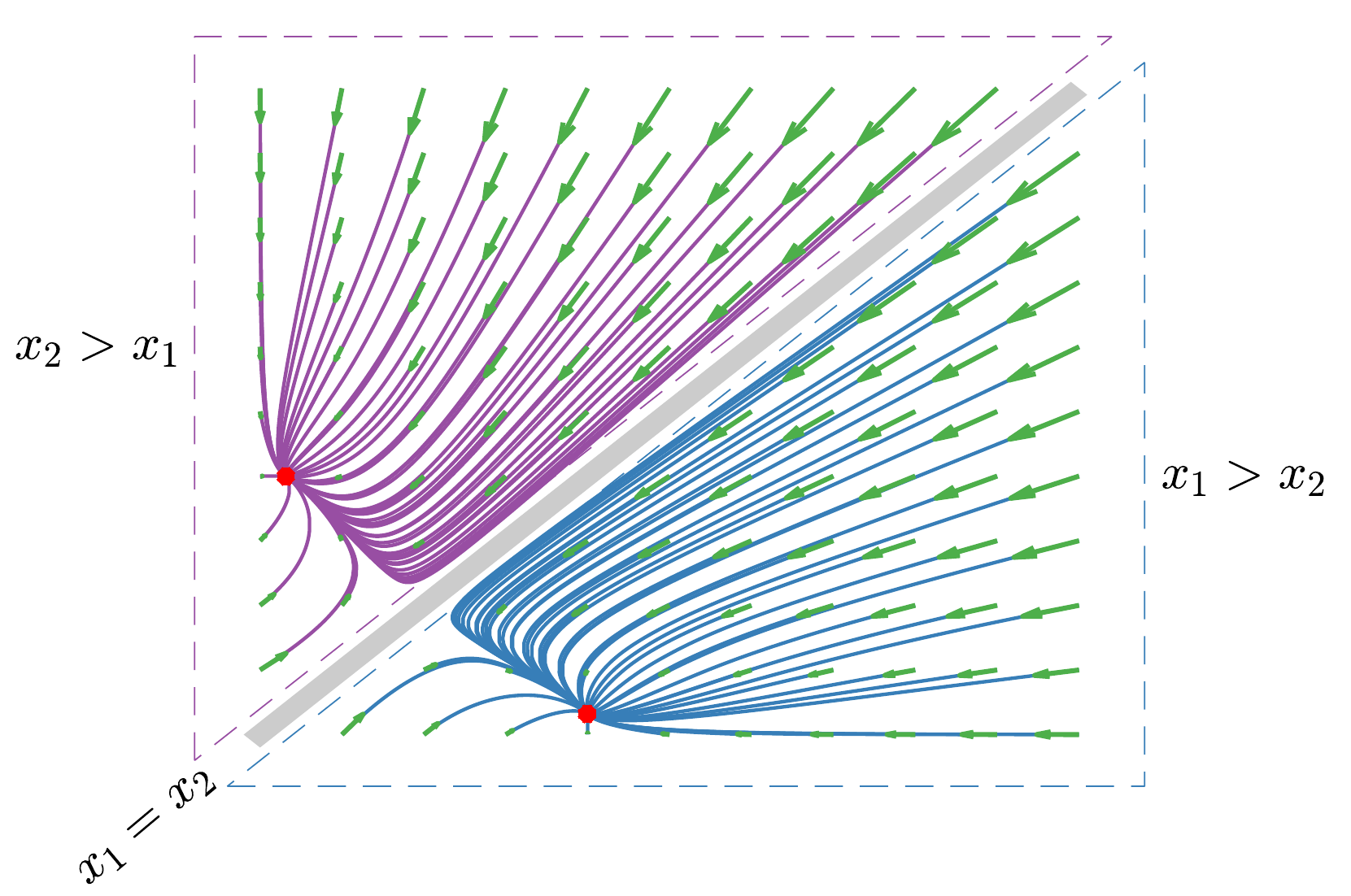}
\caption{The plane $x_1 = x_2$ acts as a mirror and the phase portrait in the region $x_2 > x_1$ is a mirror reflection of the phase portrait in the region $x_1 > x_2$.}
\label{fig:bistable_phase_portrait_symmetry}
\end{center}
\end{figure}

This system exhibits bistability, that is, this system has two stable equilibrium points and an unstable equilibrium point. The system has two invariant sets and the line $x_1=x_2$ is the separatrix that separates the two invariant sets. The phase portrait of this system is shown in  Fig. \ref{fig:bistable_phase_portrait_symmetry}.


With the given values of the parameters, the system equations are symmetric under the transformation
\begin{align}\label{reflection_sym}
(x_1,x_2)^\top\xmapsto{\g} (x_2,x_1)^\top.
\end{align}
Hence the symmetry group of the bistable toggle switch is $\mathbb{Z}_2$, where the action of the non-identity element of $\mathbb{Z}_2$ is given by Eq. (\ref{reflection_sym}). In the phase space, this corresponds to a reflection about the $x_1=x_2$ line and the $2$-dimensional representation of the non-identity element of the symmetry group is 
\begin{equation}\label{gamma_bistable}
\gamma =  \begin{pmatrix} 0 &\; 1 \\ 1 &\; 0 \end{pmatrix}.
\end{equation}

The goal is to construct the global Koopman operator using only the time-series data from any one of the invariant sets. This is achieved in three steps described below: 
\begin{enumerate}
\item Corresponding to the time-series data from any one of the invariant set, identify the dynamics applying data-driven operator theoretic methods discussed in section \ref{subsection_EDMD}. 
\item Identify the dynamics of the other invariant set by noticing that the bistable toggle switch system has reflective symmetry and applying the results from section \ref{subsection_global_koopman_equivariant}. 
\item Once the Koopman operators for each invariant set are identified, the global Koopman operator is computed by using the results of the section \ref{subsection_global_koopman_construction}. 
\end{enumerate}
To demonstrate the proposed framework, we collected time-series data from only the invariant given by the region $x_1 > x_2 $. The local Koopman operator, obtained using the DMD algorithm, is given by
\begin{align*}
\K_{right} = & \begin{pmatrix} 0.6039 & 0.0313 \\ -0.4784 & 1.0375 \end{pmatrix}.
\end{align*}

Hence, from theorem \ref{K_i_K_j_theorem}, the Koopman operator corresponding to the region $x_2>x_1$ can be identified as
\begin{align*}
\K_{left} = & \g^{-1}\K_{right}\g=\begin{pmatrix}
1.0375 & -0.4784 \\ 0.0313 & 0.6039
\end{pmatrix}, 
\end{align*}
where $\g$ is given by Eq. (\ref{gamma_bistable}). Hence the global Koopman operator is 
\[\K_{global}=\begin{pmatrix}
\K_{left} & 0\\
0 & \K_{right}
\end{pmatrix}.\]
The phase portrait corresponding to the two regions is shown in Fig. \ref{fig:bistable_phase_portrait_symmetry}. 
Moreover, as a verification, we computed the Koopman operator using data from the region $x_2>x_1$ and it was equal to 
\begin{align*}
\K_{x_2>x_1} = \begin{pmatrix}
1.0375 & -0.4784 \\ 0.0313 & 0.6039
\end{pmatrix} = \K_{left}.
\end{align*}

\subsection{Rotational Symmetry: Lorenz System}
Consider the Lorenz system as shown in Eq. \eqref{sys_lorenz}
The Lorenz system is symmetric under $\mathbb{Z}_2$ action given by 
\begin{equation}
(x,y,z)^\top\xmapsto{\g} (-x,-y,z)^\top,    
\end{equation}
which corresponds to a rotation of $180^\circ$ about the $z$-axis and matrix representation of $\g$ is
\begin{equation}\label{sym_lorenz_sim}
\begin{aligned}
\gamma = \begin{pmatrix}
-1 & 0 & 0\\
0 & -1 & 0\\
0 & 0 & 1\end{pmatrix}.
\end{aligned}
\end{equation}
\begin{figure}[h!]
\begin{center}
\includegraphics[width = 0.95 \columnwidth]{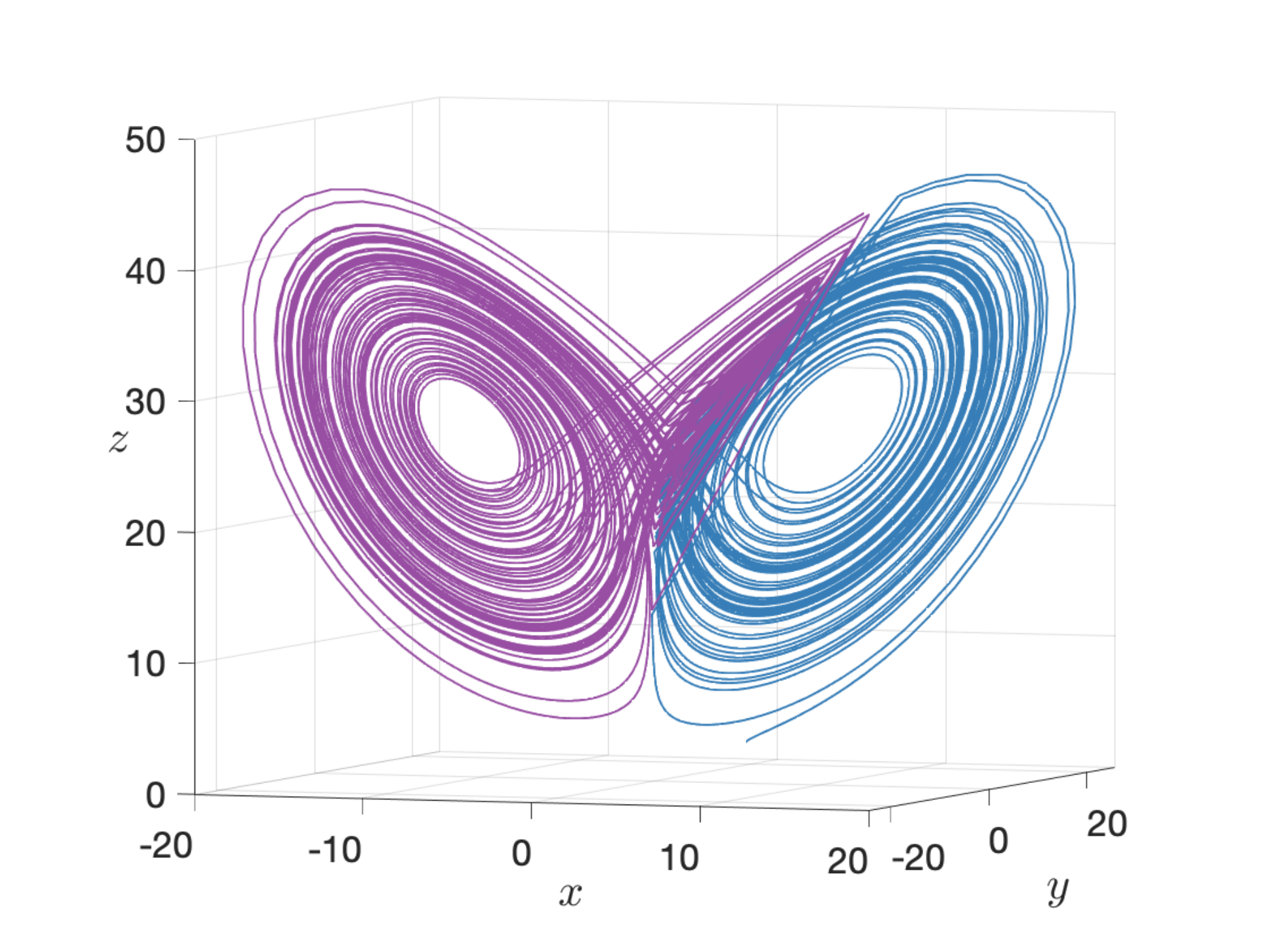}
\caption{Phase portrait of the Lorenz system.}
\label{fig:lorenz_phase_portrait_symmetry}
\end{center}
\end{figure}
The phase portrait of the Lorenz system with $\rho = 28$, $\sigma = 10$ and $\beta = 8/3$, is shown in Fig. \ref{fig:lorenz_phase_portrait_symmetry}. The colours blue and magenta correspond to the symmetric components of the strange attractor. The Koopman operator computed, using DMD algorithm, from the blue region of the attractor is
\[
\K_{blue} = \begin{pmatrix}
0.076  &  0.709  &  0.042 \\
-0.667  &  1.064  &  0.124 \\
-0.422  &  0.926  &  0.836 
\end{pmatrix}. \]
Hence, from corollary \ref{K_i_K_j_corollary}, the Koopman on the symmetric counterpart of the blue region will be
\[
\K_{magenta} = \g^{-1}\K_{blue}\g\\
=\begin{pmatrix}
0.076 &   0.709 &  -0.042 \\
-0.667 &   1.064  &  -0.124 \\ 
0.422 &  -0.926  &  0.836 \end{pmatrix}
\]
which is the same Koopman operator obtained with DMD algorithm with data points on the magenta region of the attractor.

\subsection{Reflection and Rotational Symmetry: A Hamiltonian System}
\begin{figure}[h!]
\begin{center}
\includegraphics[width = 0.95\columnwidth]{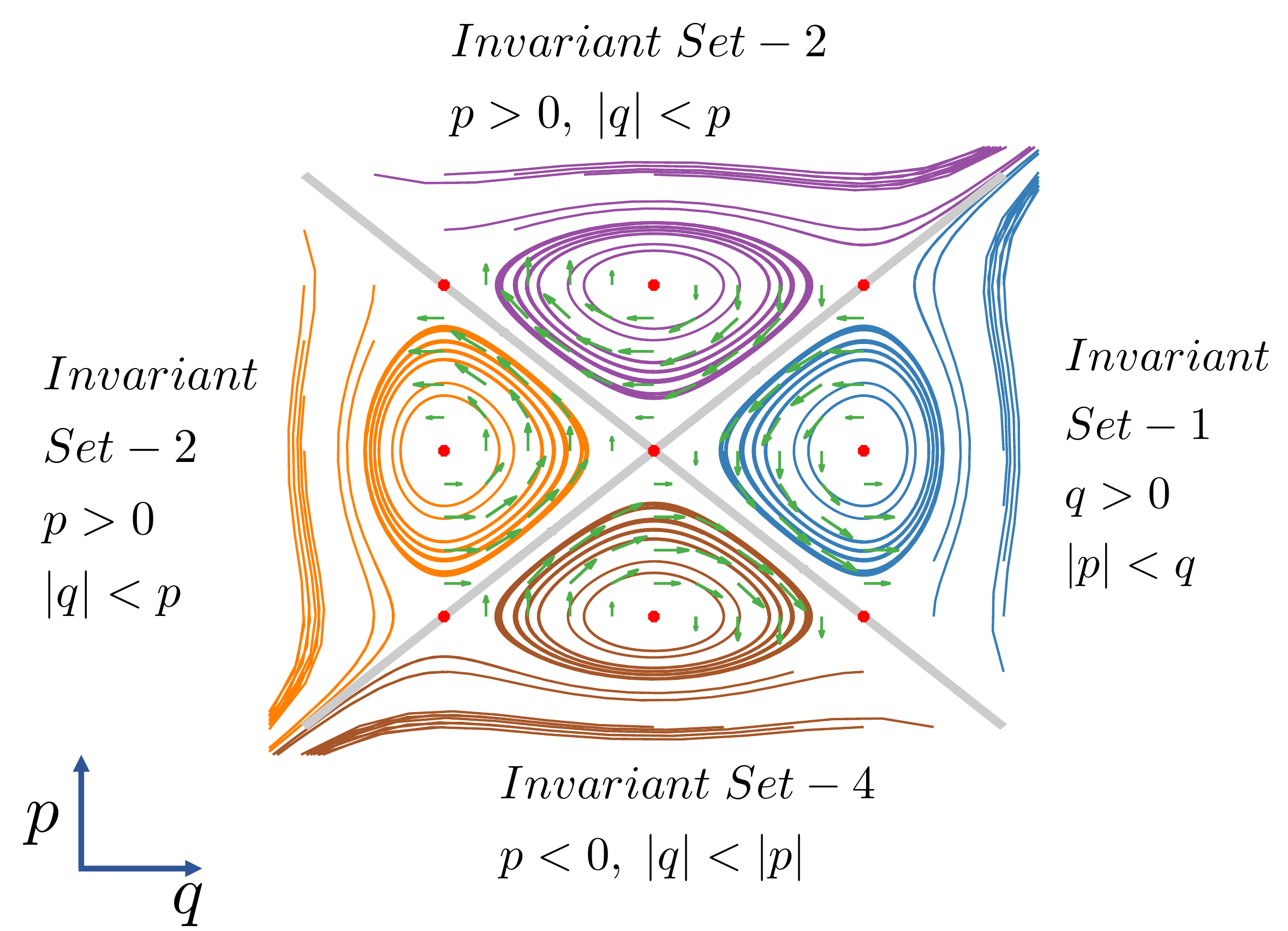}
\caption{Phase portrait of the Hamiltonian system. Suppose if we start from the invariant set - 1 (IS-1), then the invariant set - 2 (IS-2) is obtained by $\gamma_1$-action on IS-1; invariant set - 3 (IS-3) is obtained by $\gamma_2$-action on IS-1; invariant set - 4 (IS-4) is obtained by $\gamma_3$-action on IS-1, which is same as $\g_1\g_2$-action on IS-1. }
\label{fig:klein_phase_portrait_symmetry}
\end{center}
\end{figure}

Consider a Hamiltonian system with a Hamiltonian
\[H(q,p)=\frac{1}{4}p^4-\frac{9}{2}p^2 -\frac{1}{4}q^4+\frac{9}{2}q^2.\]
Hence the equations of motion are
\begin{equation}\label{klein_sym_sys}
\begin{aligned}
\dot{q} = \frac{\partial H}{\partial p} = p^3 - 9p;\quad \dot{p} = -\frac{\partial H}{\partial q} = q^3 - 9q.
\end{aligned}
\end{equation}
This system has $4$ invariant sets and the corresponding phase portrait of the system is shown in Fig. \ref{fig:klein_phase_portrait_symmetry} and the system is symmetric under the actions given by
\begin{align}
\begin{pmatrix}
q\\
p
\end{pmatrix}\xmapsto{\g_1}\begin{pmatrix}
p\\
q
\end{pmatrix},\begin{pmatrix}
q\\
p
\end{pmatrix}\xmapsto{\g_2}\begin{pmatrix}
-q\\
-p
\end{pmatrix},\begin{pmatrix}
q\\
p
\end{pmatrix}\xmapsto{\g_3}\begin{pmatrix}
-p\\
-q
\end{pmatrix}.
\end{align}

From the action of $\g_i$'s, we have
\[\g_1^2 = \g_2^2=\g_3^3=I_2, \textnormal{ and } \g_1\g_2=\g_3.\]
Hence the symmetry group of the system is the Klein 4-group $\mathbb{Z}_2\times\mathbb{Z}_2$, which has the presentation
\[\Gamma = \langle \g_1,\g_2|\g_1^2 = \g_2^2=(\g_1\g_2)^2=I_2\rangle.\]
and the matrix representation of the group elements are
\begin{align*}
\g_1 = \begin{pmatrix}
0 & 1\\
1 & 0
\end{pmatrix};\g_2 = \begin{pmatrix}
-1 & 0\\
0 & -1
\end{pmatrix};\g_3 = \begin{pmatrix}
0 & -1\\
-1 & 0
\end{pmatrix}.
\end{align*}

The local Koopman operators obtained using data from each of the invariant subspaces are
\begin{align*}
& \K_{IS-1} = \begin{pmatrix}
0.955  &  0.486\\
-0.059 &  0.215
\end{pmatrix}; 
\K_{IS-2} = \begin{pmatrix}
0.215  & -0.059\\
0.486  &  0.955
\end{pmatrix}\\
&\K_{IS-3} = \begin{pmatrix}
0.957  &  0.511\\
-0.061  &  0.214
\end{pmatrix};
\K_{IS-4} = \begin{pmatrix}
0.214  & -0.061\\
0.511  &  0.957
\end{pmatrix}
\end{align*}
and it can be seen that $\K_{IS-2}=\g_1^{-1}\K_{IS-1}\g_1$. Similar relations are found to hold true for the other local Koopman operators. Hence if only one local Koopman operator is computed from data, all the other local Koopman operators can be computed using the relation of theorem \ref{K_i_K_j_theorem}, without using data from the other invariant subspaces and they are stitched together to obtain the global Koopman operator as described in \cite{global_koopman_sai_arxiv}.

\section{Conclusions}\label{section_conclusions}
In this paper, we developed Koopman operator theoretic based methods to study the global phase space in equivariant dynamical systems. In particular, we showed that the invariant subspaces are mapped to invariant subspaces and eigenspaces are left invariant under the group action of the symmetry group and established the properties of the Koopman operator for an equivariant dynamical system under group action. Assuming the knowledge of the type of symmetry in a dynamical system, it is shown that the global phase space can be studied based on data from any one invariant subspace only.  The proposed framework is demonstrated on three different systems that possess various symmetries, such as reflective, rotational, or both. Future efforts focus on identifying the type of symmetries in a dynamical system, given the data for the global phase space.


\bibliography{subhrajit_equivariant_Koopman}            

\begin{thebibliography}{35}
\providecommand{\natexlab}[1]{#1}
\providecommand{\url}[1]{\texttt{#1}}
\providecommand{\urlprefix}{URL }
\expandafter\ifx\csname urlstyle\endcsname\relax
  \providecommand{\doi}[1]{doi:\discretionary{}{}{}#1}\else
  \providecommand{\doi}{doi:\discretionary{}{}{}\begingroup
  \urlstyle{rm}\Url}\fi

\bibitem[{Budisic et~al.(2012)Budisic, Mohr, and Mezic}]{mezic_koopmanism}
Budisic, M., Mohr, R., and Mezic, I. (2012).
\newblock Applied koopmanism.
\newblock \emph{Chaos}, 22, 047510--32.

\bibitem[{Chossat and Golubitsky(1988)}]{chossat1988symmetry}
Chossat, P. and Golubitsky, M. (1988).
\newblock Symmetry-increasing bifurcation of chaotic attractors.
\newblock \emph{Physica D: Nonlinear Phenomena}, 32(3), 423--436.

\bibitem[{Dellnitz and Junge(1999)}]{Dellnitz_Junge}
Dellnitz, M. and Junge, O. (1999).
\newblock On the approximation of complicated dynamical behavior.
\newblock \emph{SIAM Journal on Numerical Analysis}, 36, 491--515.

\bibitem[{Dellnitz et~al.(2005)Dellnitz, Junge, and et~al}]{Dellnitztransport}
Dellnitz, M., Junge, O., and et~al (2005).
\newblock Transport in dynamical astronomy and multibody problems.
\newblock \emph{International Journal of Bifurcation and Chaos}, 15, 699--727.

\bibitem[{Field(1970)}]{field1970equivariant}
Field, M. (1970).
\newblock Equivariant dynamical systems.
\newblock \emph{Bulletin of the American Mathematical Society}, 76(6),
  1314--1318.

\bibitem[{Field(1980)}]{field1980equivariant}
Field, M. (1980).
\newblock Equivariant dynamical systems.
\newblock \emph{Transactions of the American Mathematical Society}, 259(1),
  185--205.

\bibitem[{Froyland(2001)}]{froyland_extracting}
Froyland, G. (2001).
\newblock Extracting dynamical behaviour via {Markov} models.
\newblock In A.~Mees (ed.), \emph{Nonlinear Dynamics and Statistics:
  Proceedings, Newton Institute, Cambridge, 1998}, 283--324. Birkhauser.

\bibitem[{Gardner et~al.(2000)Gardner, Cantor, and
  Collins}]{gardner2000construction}
Gardner, T.S., Cantor, C.R., and Collins, J.J. (2000).
\newblock Construction of a genetic toggle switch in escherichia coli.
\newblock \emph{Nature}, 403(6767), 339.

\bibitem[{Golubitsky et~al.(2012)Golubitsky, Stewart, and
  Schaeffer}]{golubitsky2012singularities}
Golubitsky, M., Stewart, I., and Schaeffer, D.G. (2012).
\newblock \emph{Singularities and groups in bifurcation theory}, volume~2.
\newblock Springer Science \& Business Media.

\bibitem[{Horn and Johnson(2012)}]{horn2012matrix}
Horn, R.A. and Johnson, C.R. (2012).
\newblock \emph{Matrix analysis}.
\newblock Cambridge university press.

\bibitem[{Johnson and Yeung(2018)}]{johnson2018class}
Johnson, C.A. and Yeung, E. (2018).
\newblock A class of logistic functions for approximating state-inclusive
  koopman operators.
\newblock In \emph{2018 Annual American Control Conference (ACC)}, 4803--4810.
  IEEE.

\bibitem[{Junge and Osinga(2004)}]{Junge_Osinga}
Junge, O. and Osinga, H. (2004).
\newblock A set oriented approach to global optimal control.
\newblock \emph{ESAIM: Control, Optimisation and Calculus of Variations},
  10(2), 259--270.

\bibitem[{Lasota and Mackey(1994)}]{Lasota}
Lasota, A. and Mackey, M.C. (1994).
\newblock \emph{Chaos, Fractals, and Noise: Stochastic Aspects of Dynamics}.
\newblock Springer-Verlag, New York.

\bibitem[{Mauroy and Mezic(2013)}]{mezic_koopman_stability}
Mauroy, A. and Mezic, I. (2013).
\newblock A spectral operator-theoretic framework for global stability.
\newblock In \emph{Proc. of IEEE Conference of Decision and Control}. Florence,
  Italy.

\bibitem[{Mehta and Vaidya(2005)}]{Mehta_comparsion_cdc}
Mehta, P.G. and Vaidya, U. (2005).
\newblock On stochastic analysis approaches for comparing dynamical systems.
\newblock In \emph{Proceeding of IEEE Conference on Decision and Control},
  8082--8087. Spain.

\bibitem[{Mesbahi et~al.(2019)Mesbahi, Bu, and Mesbahi}]{mesbahi_symmetry}
Mesbahi, A., Bu, J., and Mesbahi, M. (2019).
\newblock On modal properties of the koopman operator for nonlinear systems
  with symmetry.
\newblock In \emph{2019 American Control Conference (ACC)}, 1918--1923. IEEE.

\bibitem[{Mezic and Banaszuk(2000)}]{Mezic2000}
Mezic, I. and Banaszuk, A. (2000).
\newblock Comparison of systems with complex behavior: spectral methods.
\newblock In \emph{Proceedings of the 39th IEEE Conference on Decision and
  Control}, 1224--1231.

\bibitem[{Mezi\'{c} and Banaszuk(2004)}]{Mezic_comparison}
Mezi\'{c}, I. and Banaszuk, A. (2004).
\newblock Comparison of systems with complex behavior.
\newblock \emph{Physica D}, 197, 101--133.

\bibitem[{Mezi{\'c}(2005)}]{mezic2005spectral}
Mezi{\'c}, I. (2005).
\newblock Spectral properties of dynamical systems, model reduction and
  decompositions.
\newblock \emph{Nonlinear Dynamics}, 41(1-3), 309--325.

\bibitem[{Nandanoori et~al.(2019)Nandanoori, Sinha, and
  Yeung}]{global_koopman_sai_arxiv}
Nandanoori, S.P., Sinha, S., and Yeung, E. (2019).
\newblock Data-driven operator theoretic methods for global phase space
  learning.
\newblock \emph{arXiv preprint arXiv:1910.03011}.

\bibitem[{Petersen(1989)}]{petersen_ergodic_book}
Petersen, K.E. (1989).
\newblock \emph{Ergodic theory}, volume~2.
\newblock Cambridge University Press.

\bibitem[{Raghunathan and Vaidya(2014)}]{raghunathan2014optimal}
Raghunathan, A. and Vaidya, U. (2014).
\newblock Optimal stabilization using lyapunov measures.
\newblock \emph{IEEE Transactions on Automatic Control}, 59(5), 1316--1321.

\bibitem[{Salova et~al.(2019)Salova, Emenheiser, Rupe, Crutchfield, and
  D’Souza}]{koopman_symmetry}
Salova, A., Emenheiser, J., Rupe, A., Crutchfield, J.P., and D’Souza, R.M.
  (2019).
\newblock Koopman operator and its approximations for systems with symmetries.
\newblock \emph{Chaos: An Interdisciplinary Journal of Nonlinear Science},
  29(9), 093128.

\bibitem[{Schmid(2010)}]{DMD_schmitt}
Schmid, P.J. (2010).
\newblock Dynamic mode decomposition of numerical and experimental data.
\newblock \emph{Journal of Fluid Mechanics}, 656, 5--28.

\bibitem[{Sinha et~al.(2018{\natexlab{a}})Sinha, Bowen, and
  Vaidya}]{robust_DMD_arxiv}
Sinha, S., Bowen, H., and Vaidya, U. (2018{\natexlab{a}}).
\newblock On robust computation of koopman operator and prediction in random
  dynamical systems.
\newblock \emph{arXiv preprint arXiv:1803.08562}.

\bibitem[{Sinha et~al.(2018{\natexlab{b}})Sinha, Huang, and
  Vaidya}]{robust_DMD_ACC}
Sinha, S., Huang, B., and Vaidya, U. (2018{\natexlab{b}}).
\newblock Robust approximation of koopman operator and prediction in random
  dynamical systems.
\newblock In \emph{2018 Annual American Control Conference (ACC)}, 5491--5496.
  IEEE.

\bibitem[{Sinha et~al.(2019{\natexlab{a}})Sinha, Nandanoori, and
  Yeung}]{sinha_online_koopman_arxiv}
Sinha, S., Nandanoori, S.P., and Yeung, E. (2019{\natexlab{a}}).
\newblock Online learning of dynamical systems: An operator theoretic approach.
\newblock \emph{arXiv preprint arXiv:1909.12520}.

\bibitem[{Sinha et~al.(2019{\natexlab{b}})Sinha, Vaidya, and
  Yeung}]{sparse_Koopman_acc}
Sinha, S., Vaidya, U., and Yeung, E. (2019{\natexlab{b}}).
\newblock On computation of koopman operator from sparse data.
\newblock In \emph{2019 American Control Conference (ACC)}, 5519--5524. IEEE.

\bibitem[{Sparrow(2012)}]{sparrow2012lorenz}
Sparrow, C. (2012).
\newblock \emph{The Lorenz equations: bifurcations, chaos, and strange
  attractors}, volume~41.
\newblock Springer Science \& Business Media.

\bibitem[{Surana and Banaszuk(2016)}]{surana_observer}
Surana, A. and Banaszuk, A. (2016).
\newblock Linear observer synthesis for nonlinear systems using koopman
  operator framework.
\newblock \emph{IFAC-PapersOnLine}, 49(18), 716--723.

\bibitem[{Susuki and Mezic(2011)}]{susuki2011nonlinear}
Susuki, Y. and Mezic, I. (2011).
\newblock Nonlinear koopman modes and coherency identification of coupled swing
  dynamics.
\newblock \emph{IEEE Transactions on Power Systems}, 26(4), 1894--1904.

\bibitem[{Vaidya and Mehta(2008)}]{Vaidya_TAC}
Vaidya, U. and Mehta, P.G. (2008).
\newblock Lyapunov measure for almost everywhere stability.
\newblock \emph{IEEE Transactions on Automatic Control}, 53(1), 307--323.

\bibitem[{Williams et~al.(2015)Williams, Kevrekidis, and
  Rowley}]{EDMD_williams}
Williams, M.O., Kevrekidis, I.G., and Rowley, C.W. (2015).
\newblock A data--driven approximation of the koopman operator: Extending
  dynamic mode decomposition.
\newblock \emph{Journal of Nonlinear Science}, 25(6), 1307--1346.

\bibitem[{Yeung et~al.(2017)Yeung, Kundu, and Hodas}]{yeung2017learning}
Yeung, E., Kundu, S., and Hodas, N. (2017).
\newblock Learning deep neural network representations for koopman operators of
  nonlinear dynamical systems.
\newblock \emph{arXiv preprint arXiv:1708.06850}.

\bibitem[{Yeung et~al.(2018)Yeung, Liu, and Hodas}]{yeung2018koopman}
Yeung, E., Liu, Z., and Hodas, N.O. (2018).
\newblock A koopman operator approach for computing and balancing gramians for
  discrete time nonlinear systems.
\newblock In \emph{2018 Annual American Control Conference (ACC)}, 337--344.
  IEEE.

\end{thebibliography}

\end{document}